               \documentclass[11pt,onecolumn]{article}                                           

\usepackage{algorithm, algorithmic, setspace}
\usepackage{graphicx} 
\usepackage{amsmath} 
\usepackage{amssymb}  
\usepackage{amsthm}
\usepackage{amsfonts}
\usepackage{color}
\usepackage{subfig}
\usepackage[english]{babel}
\usepackage[latin1]{inputenc}

\newcommand{\fun}{\mathcal{F}}
\newcommand{\lag}{\mathcal{L}}

\newcommand{\soft}{\mathbb{S}}

\newcommand{\Sc}{\bar{S}}

\newcommand{\D}{[-d,d]^n}
\newcommand{\cone}{\mathcal{C}}

\newcommand{\argmin}[1]{\underset{#1}{\mathrm{argmin\,}}}

\newcommand{\xmin}{x^{\star}}
\newcommand{\xgmin}{x^{\star\star}}

\newcommand{\la}{\lambda}
\newcommand{\g}{\mathrm{g}}
\newcommand{\si}{\mathrm{sign}}
\newcommand{\xtrue}{\widetilde{x}}
\newcommand{\noise}{\eta}

\newcommand{\ball}{\mathcal{B}_{\varepsilon}(\xmin)}

\newcommand{\R}{\mathbb{R}}

\newtheorem{lemma}{Lemma}
\newtheorem{proposition}{Proposition}
\newtheorem{theorem}{Theorem}
\newtheorem{corollary}{Corollary}
\newtheorem{remark}{Remark}

\newtheorem{definition}{Definition}
\newtheorem{assumption}{Assumption}
\graphicspath{{figures/}}

\title{Sparse learning with concave regularization:\\ relaxation of the irrepresentable condition}

\author{V. Cerone, S. M. Fosson$^{*}$, D. Regruto, A. Salam \thanks{$^{*}$ Corresponding author. The authors are with the Dipartimento di Automatica e Informatica,
Politecnico di Torino, corso Duca degli Abruzzi 24, 10129 Torino, Italy;
    e-mail: sophie.fosson@polito.it.
}
}
\date{}
\begin{document}
\maketitle
\thispagestyle{empty}
\pagestyle{empty}

\begin{abstract}
Learning sparse models from data is an important task in all those frameworks where relevant information should be identified within a large dataset. This can be achieved by formulating and solving suitable sparsity promoting optimization problems. As to linear regression models,  Lasso is the most popular convex approach, based on an $\ell_1$-norm regularization. In contrast,  in this paper, we analyse a concave regularized approach, and we prove that it relaxes the irrepresentable condition, which is sufficient and essentially necessary for Lasso to select the right significant parameters. In practice, this has the benefit of reducing the number of necessary measurements with respect to Lasso. Since the proposed problem is non-convex, we also discuss different algorithms to solve it, and we illustrate the obtained enhancement via numerical experiments.
\end{abstract}
\section{Introduction}\label{sec:in}
Sparse learning is the science of building parsimonious models from data. The main motivation for sparse learning is the concrete need of extracting relevant information from large collections of data, which nowadays are commonly available in many scientific fields. This task prevents drawbacks such as overfitting, redundancies, numerical complexity, and scarce understanding of the physical behavior of systems; we refer the reader to \cite{has09book,has15book,bru19book} for a comprehensive illustration of these issues. Recent applications of sparse learning can be found in identification of linear and non-linear systems, see \cite{car13,fat18} and \cite{bru16,bru19book}, respectively; in model predictive control, see \cite{gal16book}; in neural networks and deep learning, see \cite{tar18,bol19}. In signal processing, the exploitation of sparsity has a long history, provided that many signals (e.g., images) admit sparse representations in opportune bases. In this context, the compressed sensing (CS) theory has been developed in the last fifteen years, which states that a sparse vector $x\in\R^n$ can be recovered from compressed, possibly noisy, linear measurements, see \cite{don06,fou13}. Beyond signal processing, the CS paradigm has been applied, e.g., to linear/non-linear system identification, see \cite{tot11,san11,car13, bru19book}.
 
Often, sparse learning is formulated as an optimization problem where sparsity is promoted by a suitable regularization term. In principle, the $\ell_0$-norm, i.e., the number of non-zero components of a vector, is the correct pseudonorm to represent the sparsity level; nevertheless, the $\ell_0$-norm is non-convex and leads to NP-hard optimization. As an alternative, the $\ell_1$-norm has been studied and proven to be the best convex approximation of the $\ell_0$-norm. In the linear regression setting, the Lasso problem, introduced in \cite{tib96}, is very popular, and consists in the minimization of an $\ell_1$ regularized, least squares cost function.

The reliability of the Lasso estimator has been largely studied in the literature, in particular for  what concerns its variable-selection consistency (VSC), i.e., the ability of identifying the support, which is the set of non-zero/significant components of the unknown vector. Finding the support is the most important task in sparse learning and CS, since its knowledge is sufficient to recover the complete vector. As illustrated in \cite{fuc04,zha06,wai09}, Lasso enjoys the VSC if the so-called irrepresentable condition holds. More precisely, as we illustrate in the next section, this condition is sufficient and essentially necessary, see \cite{bul11}.

In this paper, we analyse an alternative approach to Lasso for sparse learning, based on a concave, semi-algebraic regularization. Recently, the use of non-convex regularizers has been gaining attention in the literature, see, e.g., \cite{can08rew,woo16,sel17,fox19cdc}. The rationale is that the shape of a non-convex regularizer is closer to the $\ell_0$-norm than the $\ell_1$-norm. Numerical experiments show that the non-convex approach is usually more effective than Lasso. Nevertheless, theoretical results are still missing, in particular in the compressed framework.

The goal of this paper is to prove that the considered non-convex estimator is more effective than Lasso, in the sense that it enjoys the VSC under a relaxed irrepresentable condition, even in a CS setting. This result is obtained by exploiting a restricted eigenvalue property and a boundedness assumption. Under these hypotheses, a minimum of the proposed functional has the correct support; sufficient conditions to evaluate whether this minimum is local or global are provided. In practice, the relaxation of the irrepresentable condition implies a reduction of the number of necessary measurements to select the right parameters.  
Since the problem is semi-algebraic, the Lasserre's approach can be used to compute the minimum \cite{lasbook}. On the other hand, iterative algorithms can be used for local minimization and are effective in many cases. This is analysed through numerical simulations.

The paper is organized as follows. In Sec. \ref{sec:ps}, we present the problem and we illustrate the considered cost function. In Sec. \ref{sec:bac}, we review how to derive the irrepresentable condition for Lasso. This background opens the way to understanding the analysis of the proposed approach, developed in Sec. \ref{sec:ta}, which is the core of the contribution. In Sec. \ref{sec:nr}, we discuss the minimization algorithms to implement the proposed approach and we show some numerical results. Finally, we draw some conclusions in Sec. \ref{sec:con}.
\section{Problem statement and background}\label{sec:ps}
In this paper, we consider the following sparse optimization problem. We aim to estimate a $k$-sparse vector $\xtrue\in\R^n$ from noisy linear measurements:
\begin{equation}\label{theproblem}
 y=A\xtrue+\noise
\end{equation}
where $A\in\R^{m,n}$, and $\eta\in\R^m$ is a measurements noise. In particular, we focus on the CS case $m<n$, which is more challenging. We call $S\subset \{1,\dots,n\}$ the support of $\xtrue$, and $\Sc$ its complementary. The cardinality of $S$ is $|S|=k$. It is natural to assume that the number of measurements is not smaller than $|S|$, that is, $k\leq m$.

For any $v\in\R^n$, we denote by $v_S\in\R^k$ the vector that contains the components of $v$ restricted to $S$. Similarly, $A_S\in\R^{m,k}$ is the submatrix of $A$ formed by the columns of $A$ indexed in $S$.

We consider the following assumptions.
\begin{assumption}\label{ass:bounded}
The non-zero parameters are bounded, i.e., $\xtrue_i\in[-d,d]$ for each $i \in S$, where $d>0$ is known.
\end{assumption}
\begin{assumption}\label{ass:AS}
$A_S$ has full rank.
\end{assumption}
Assumption \ref{ass:bounded} is realistic because in real-world applications a prior estimate on the maximum magnitude is usually available. By considering $k\leq m$, Assumption \ref{ass:AS} states that the columns of $A_S$ are linearly independent, which is a standard condition in CS; see, e.g., \cite{tib13,fou13,has15book} for details.
Due to noise and compression, the general approach to the considered problem is a regularized least squares. In particular,  Lasso \cite{tib96} reads as follows:
\begin{equation}\label{lasso}  
{\text{\bf{Lasso}}}:~~~ \min_{x\in\D}\frac{1}{2}\left\|y-A x \right\|_2^2+\la\|x\|_1
\end{equation}
where $\la>0$ is a design parameter which can be assessed based on prior information on $k$, see \cite{fou13}.
In the literature, the VSC of Lasso is analysed by using the irrepresentable condition.

\begin{definition}\label{def:irr}{\bf{Irrepresentable condition (IRR)}}
We say that $A$ satisfies the $(\omega,S)$ irrepresentable condition,  $(\omega,S)$-IRR for short, if there exists $\omega>0$ such that
\begin{equation}\label{irr}
\left\|A_S^{\dagger} A_{\Sc}\right\|_{1}< \frac{1}{\omega}
\end{equation}
 \end{definition}
where
\begin{equation}
 A_S^{\dagger}=(A_S^\top A_S)^{-1}A_S^\top\in\R^{k,m}
\end{equation}
is the Moore-Penrose pseudoinverse of $A_S$; we notice that  $A_S^\top A_S$ is invertible by Assumption \ref{ass:AS}. 

Slightly different irrepresentable conditions are introduced and analyzed in the literature; see, e.g., \cite{fuc04,fuc05,mei06,zha06,wai09}. We refer the reader to Sec. \ref{sec:bac} for a discussion. The key common idea is that  $(\omega,S)$-IRR with $\omega>1$, also known as strong irrepresentable condition, is sufficient for the VSC of Lasso if some assumptions on the measurement noise hold and provided a suitable tuning of $\la$. On the other hand, $(\omega,S)$-IRR with $\omega=1$, also known as weak irrepresentable condition, is necessary for VSC. This is studied both for the asymptotic case $m\to\infty$, see \cite{mei06,zha06}, and non-asymptotic case, see \cite{fuc04,fuc05,wai09}.

%
%

In this paper, we consider a concave alternative to $\ell_1$ regularization. In the last years, different concave penalties have been proposed in sparse optimization, in particular, $\log(x)$, see \cite{can08rew}; $\ell_p$, with $p\in(0,1)$, see \cite{fou09}, and minimax concave penalties (MCP, \cite{zha10MCP}). In some cases, the use of such regularizers is associated with the exploitation of $\ell_1$-reweighting techniques for minimization, see \cite{can08rew}. The key idea behind the use of concave penalties is that they are closer to the $\ell_0$-norm with respect to the $\ell_1$-norm, therefore they are expected to promote sparsity more accurately. We refer the reader to \cite{woo16,fox18} for their use in CS.

In this paper, we propose an MCP-based shrinkage and selection method, denoted as MCPS$^2$, which reads as follows.
\begin{equation}\label{mcps2}  
\begin{split}
\text{\bf{MCPS}}^2:~~~ &\min_{x\in\D}\fun(x):=\frac{1}{2}\left\|y-A x \right\|_2^2+\la r(x)\\
 {\text{ where }}~~~&r(x):=d\left\|x\right\|_1-\frac{1}{2}\left\|x \right\|_2^2,~~\la>0.
\end{split}
\end{equation}
MCPS$^2$ is partially studied in \cite{fox18tsp}, limited to CS of finite-valued signals, while in  \cite{fox19cdc} an MCP cost function is tested for the recovery of non-negative signals from low-precision data. Differently from those works, in this paper we extend the analysis to any real signal and we prove novel theoretical results.

For our analysis, we exploit the restricted eigenvalue condition. Let us define the cone
\begin{equation}\label{cone}
\cone(\alpha,S):=\{x\in\R^n: \left\|x_{\Sc}\right\|_1\leq \alpha \left\|x_{S}\right\|_1\}.
\end{equation}

\begin{definition}\label{def:re}{\bf{Restricted eigenvalue condition (RE)}}, see \cite[Section 11.2.2]{has15book}).
We say that $A$ satisfies $(\alpha,\phi,S)$-RE  if there exist $\alpha>0$ and $\phi>0$ such that
 \begin{equation}
\frac{\|Av\|^2_2}{\|v\|_2^2}\geq\phi ~{\text{for any}}~ v\in\cone(\alpha,S)\setminus\{0\}.
\end{equation}
\end{definition}

We remark that RE is generally weaker than IRR and than the restricted isometry property (RIP, see, e.g., \cite[Definition 10.2]{has15book}), the latter being largely used in CS; we refer the reader to \cite{van09} for a thorough study of the relationships between the different  conditions considered in CS and sparse optimization. Moreover, while  RIP is proven to hold for random, independent matrices,  RE is proven to hold for a wider class of random, correlated matrices, see \cite{ras10}. Thus, RE matches with a larger number of applications, such as autoregressive models.

We also remark that  RE is used to evaluate the $\ell_2$ error of Lasso (but not its VSC): Theorem 11.1 in \cite{has15book} states that if $(3,\phi,S)$-RE holds and $\|\noise\|_{\infty}\leq m\la$, then $\|z-\xtrue\|_2\leq \frac{3\omega}{\omega-1}\sqrt{k}\la$, where $z$ is the Lasso solution. As explained in \cite[Section 11.2.2]{has15book},  RE  originates from the observation that $\|Ax-y\|_2^2$ is not strongly convex in CS, since $A^\top A$ has not full rank, and the corresponding quadratic form is positive semi-definite. For this motivation, one looks for a strong convexity restricted to a significant subset, which for Lasso is the cone $\cone(3,S)$.
\section{Irrepresentable condition for Lasso}\label{sec:bac}
Before introducing the proposed analysis, we review the role of the irrepresentable condition for the VSC of Lasso. The results that we report in this section can be retrieved in \cite{fuc04,fuc05}.

We say that $\g\in\R^n$ is a subgradient of $f:\R^n\to \R$ at $x\in\R^n$ if for all $z\in\R^n$
\begin{equation}\label{def:sub}
 f(z)\geq f(x)+\g^\top(z-x).
\end{equation}
The set of the subgradients of $f$ at $x$ is known as subdifferential $f$ at $x$ and is denoted $\partial f(x)$.
If $f$ is convex and differentiable, then its gradient is a subgradient.


The subdifferential of the $\ell_1$ norm is defined as
\begin{equation}\label{sub_abs}
\partial \|v\|_1=\{\g\in\R^n:~ \g_i=\si(v_i)~\text{ if } v_i\neq 0,~\text{ and }~ |\g_i|\leq 1~\text{ if } v_i=0\}.
\end{equation}

Since Lasso is convex and subdifferentiable, we can derive the optimality condition by evaluating the subgradient equation
\begin{equation}\label{sub_eq}
\exists \g\in\partial \|x\|_1 ~\text{ such that }~A^\top (Ax-y)+\la\g =0
\end{equation}
By splitting \eqref{sub_eq} into zero and non-zero components, we obtain the conditions
\begin{equation}\label{sub_split}
\left\{
 \begin{array}{ll}
  S\text{-condition:}~~A_i^\top (Ax-y)+\la\si(x_i) =0& \text{ if } x_i \neq 0;\\
  \Sc\text{-condition:}~~|A_i^\top (Ax-y)|\leq \la& \text{ if } x_i =0.\\
 \end{array}\right.
\end{equation}
If $A_S$ has full rank and the inequality of \eqref{sub_split} is strict, then $x$ is the unique minimizer; see \cite{fuc04,fuc05}.

The irrepresentable condition for the VSC of Lasso derives from evaluating the solutions of \eqref{sub_split} that have support $S$, that is the same support of the true $\xtrue$.

Let $\xmin\in\R^n$ be a vector with support $S$ as $\xtrue$. We wonder under which conditions $\xmin$ is a solution to \eqref{sub_split} and thus the minimizer of Lasso. We assume that $A_S$ has full rank.

From the $S$-condition in \eqref{sub_split},
\begin{equation}
\begin{split}
  &A_S^\top (A\xmin-y)+\la\si(\xmin_S) =0\\
  &\Leftrightarrow A_S^\top (A(\xmin-\xtrue)-\noise)+\la\si(\xmin_S) =0\\
  &\Leftrightarrow A_S^\top (A_S(\xmin_S-\xtrue_S)-\noise)+\la\si(\xmin_S) =0\\
  &\Leftrightarrow \xmin_S-\xtrue_S= (A_S^\top A_S)^{-1}(A_S^\top \noise-\la\si(\xmin_S)).
\end{split}
\end{equation}
From the last expression, we notice that $\si(\xmin)=\si(\xtrue)$ if the noise $\noise$ and the design hyperparameter $\la$ are sufficiently small.

More precisely,  $\|\xmin_S-\xtrue_S\|_{\infty}\leq \min_{i\in S}|\xtrue_i|$ is a sufficient condition for $\si(\xmin)=\si(\xtrue)$. Therefore, by straightforward computations,
\begin{equation}\label{lasso_limitazione_su_lambda}
 \|(A_S^\top A_S)^{-1}\|_{\infty}(\|A_S^\top \noise\|_{\infty}+\la)\leq \min_{i\in S}|\xtrue_i| ~\Rightarrow~\si(\xmin)=\si(\xtrue).
\end{equation}

Now, we replace $\xmin_S-\xtrue_S= (A_S^\top A_S)^{-1}(A_S^\top \noise-\la\si(\xmin_S))$ in the $\Sc$-condition in \eqref{sub_split}:
\begin{equation}\label{calcoli}
\begin{split}
  &\|A_{\Sc}^\top (A\xmin-y)\|_{\infty}\leq \la\\
  &\Leftrightarrow \|A_{\Sc}^\top (A_S(\xmin-\xtrue)-\noise)\|_{\infty}\leq \la\\
   &\Leftrightarrow \|A_{\Sc}^\top A_S(A_S^\top A_S)^{-1}(A_S^\top \noise-\la\si(\xmin_S))-A_{\Sc}^\top\noise\|_{\infty}\leq \la\\
    &\Leftrightarrow \|A_{\Sc}^\top A_S^{\dagger,\top}(A_S^\top \noise-\la\si(\xmin_S))-A_{\Sc}^\top\noise\|_{\infty}\leq \la\\
      &\Leftarrow \la\|A_{\Sc}^\top A_S^{\dagger\top}\si(\xmin_S)\|_{\infty}+\|A_{\Sc}^\top(A_S^{\dagger\top}A_S^\top-I)\noise\|_{\infty}\leq \la\\
\end{split}
\end{equation}
We notice that $A_S^{\dagger\top}A_S^\top=A_S A_S^{\dagger}$ and $\|A_{\Sc}^\top A_S^{\dagger\top}\|_{\infty}=\| A_S^{\dagger}A_{\Sc}\|_{1}$ By calling $\zeta(\noise):=A_{\Sc}^\top(A_S A_S^{\dagger}-I)\noise$,
\begin{equation}
\begin{split}
      &\la\|A_{\Sc}^\top A_S^{\dagger\top}\si(\xmin_S)\|_{\infty}+\|\zeta(\noise)\|_{\infty}\leq \la\\
        &\Leftarrow\la\|A_{\Sc}^\top A_S^{\dagger\top}\|_{\infty}\|\si(\xmin_S)\|_{\infty}+\|\zeta(\noise)\|_{\infty}\leq \la\\
          &\Leftrightarrow \| A_S^{\dagger}A_{\Sc}\|_{1}+\frac{\|\zeta(\noise)\|_{\infty}}{\la}\leq 1\\
\end{split}
\end{equation}
In conclusion, $\xmin$ is a minimizer of Lasso if
\begin{equation}\label{lasso_irr}
\| A_S^{\dagger}A_{\Sc}\|_{1}+\frac{\|\zeta(\noise)\|_{\infty}}{\la}\leq 1
\end{equation}
and it is the unique minimizer if the inequality is strict. In particular, if $\noise=0$, thus $\zeta(\noise)=0$, then $(\omega,S)$-IRR with $\omega>1$ is sufficient to have the minimizer $\xmin$. Since $\xmin$ has the same support of $\xtrue$, under this condition, Lasso has the VSC property.

On the other hand, from \eqref{calcoli} we see that even in the noise free case  $\|A_{\Sc}^\top A_S^{\dagger\top}\|_{\infty}\leq 1$ is necessary to have VSC for all $\xtrue\in\R^n$ with support $S$, that is, independently from the signs of the non-zero components.

\section{Relaxed irrepresentable condition for MCPS$^2$}\label{sec:ta}
In this section, we analyse the VSC of MCPS$^2$ defined in \eqref{mcps2} in a CS setting ($m<n$). We prove that a vector $\xmin$ with the same support of $\xtrue$  is a minimizer of $\fun$ in \eqref{mcps2}, under a relaxed $(\omega,S)$-IRR if compared to Lasso, i.e., $\omega$ may be smaller than 1. Therefore, a descent algorithm that starts sufficiently close to $\xmin$ can recover the right support. Furthermore, under additive conditions, $\xmin$ is proven to be the global minimum.

\subsection{Local minimum}\label{sec:ta1}
Under Assumption \ref{ass:AS},  $A_S^\top A_S$ is positive definite. Moreover, se set $\la$ so that $A_S^\top A_S-\la I$ is invertibile.

Let use denote
\begin{equation}
\mu:=\frac{1}{d}\min_{i\in S}|\xtrue_i|\in (0,1].
\end{equation}

We notice that the subdifferential arguments used in Sec. \ref{sec:bac} for Lasso cannot be used for MCPS$^{2}$ because $\fun$ in \eqref{mcps2} is nonconvex and in general the subgradient does not exist.


We define a candidate minimizer $\xmin\in\D$ as follows. First, we consider \eqref{mcps2} restricted to the true support $S$. Since $A_S^\top A_S-\la I$ is positive definite, the restricted problem is convex and we compute a minimizer by solving the subgradient equation
\begin{equation}\label{restricted_sub}
A_S^{\top}(A_S x_s-y)-\la x_s +\la d \g_S =0.
\end{equation}
Then, we set to zero the components of $\xmin$ in $\Sc$.
Thus, we have
\begin{equation}\label{defmin}
\begin{split}
&\xmin_S:=\xtrue_S+(A_S^\top A_S-\la I)^{-1}\left[\la \xtrue_S-\la d\si(\xmin_S)+A_S^\top\noise\right],\\
&\xmin_{\Sc}:=0.
\end{split}
\end{equation}
In \eqref{defmin}, we replace $\g_S$ by $\si(\xmin_S)$ because, according to the following lemma, $\si(\xtrue_S)=\si(\xmin_S)$ under some conditions. Therefore, $\xmin_i \neq 0$ for $i\in S$.

\begin{lemma}\label{the_signs}
If
\begin{equation}\label{limit_lambda}\|(A_S^\top A_S-\la I)^{-1}\|_{\infty}\left[ \la d (1-\mu)+\|A_S^\top\noise\|_{\infty}\right]\leq \mu d\end{equation}
then  $$\si(\xtrue_S)=\si(\xmin_S).$$
\end{lemma}
\begin{proof}
A sufficient condition for  $\si(\xtrue_S)=\si(\xmin_S)$ is
$$\|\xtrue_S-\xmin_S\|_{\infty}\leq \min_{i\in S}|\xtrue_i| = \mu d.$$

From \eqref{defmin},
\begin{equation*}
 \begin{split}
  \|\xtrue_S-\xmin_S\|_{\infty} &\leq \|(A_S^\top A_S-\la I)^{-1}\|_{\infty}\|\la \xtrue_S-\la d\si(\xtrue_S)+A_S^\top\noise\|_{\infty}\\
  &\leq \|(A_S^\top A_S-\la I)^{-1}\|_{\infty}\left[ \la\| \xtrue_S- d\si(\xtrue_S)\|_{\infty}+\|A_S^\top\noise\|_{\infty}\right]\\
  &\leq \|(A_S^\top A_S-\la I)^{-1}\|_{\infty}\left[ \la d (1-\mu)+\|A_S^\top\noise\|_{\infty}\right].
 \end{split}
\end{equation*}
The thesis follows by upperbounding the last expression with $\mu d$.
\end{proof}
Lemma \ref{the_signs} states that the components of $\xmin$ have the same signs of $\xtrue$ if the noise and $\la$ are sufficiently small and $\mu$ is sufficiently large. 

In the following result, we prove that $\xmin$ defined in \eqref{defmin} is a local minimizer of MCPS$^{2}$.

\begin{theorem}\label{theo1}
Let $y=A\xtrue+\noise$, where $\xtrue\in\D$ has support $S$, with $|S|=k$.

Let \eqref{limit_lambda} hold.


%

If the following conditions hold:

\begin{itemize}
 \item[C1.]
 $A$ satisfies $(\omega,S)$-IRR,
\item[C2.] $A$ satisfies $(\alpha,\phi,S)$-RE with $\phi>\la$,
\item[C3.]  $\la d-\frac{1}{\omega}\left\|A_S^\top A_S(A_S^\top A_S-\la I)^{-1}\right\|_{\infty}\left[\la d(1-\mu)+\|A_S^\top\noise\|_{\infty}\right]-\left\|A_{\Sc}^\top\noise\right\|_{\infty}>\la\varepsilon\frac{1+\alpha}{\alpha}$ where $\epsilon>0$ is arbitrarily small,
\end{itemize}
then, $\xmin$ is a local minimum of $\fun(x)$, $x\in\D$.
\end{theorem}
\begin{proof}

To prove that $\xmin$ is a local minimum, we show that
$\fun(\xmin+h)>\fun(\xmin)$ for a perturbation $h\in\R^n\setminus\{0\}$. In particular, we assume that
$\|h\|_{\infty}\leq \varepsilon$
for an arbitrarily small $\varepsilon>0$.

From \eqref{mcps2},
\begin{equation}\label{eq:start}
\begin{split}
\fun(\xmin+h)-\fun(\xmin)=&\frac{1}{2}\left\|Ah\right\|_2^2+ h^\top A^\top(A\xmin-y)+\la \left[r(\xmin+h)- r(\xmin)\right]
\end{split}
\end{equation}
where
%

\begin{equation}
\begin{split}
 r(\xmin+h)-r(\xmin)&=d\left\|\xmin+h\right\|_1- d\left\|\xmin\right\|_1-\frac{\left\|h\right\|_2^2}{2}- h^\top\xmin\\
 &=d\left\|\xmin_{S}+h_{S}\right\|_1+ d\left\|h_{\Sc}\right\|_1- d\left\|\xmin_{S}\right\|_1-\frac{\left\|h\right\|_2^2}{2}-h_S^\top\xmin_{S}.
\end{split}
\end{equation}

For each $i\in S$, we can assume that $|h_i|<|\xmin_i|$, so that $\si(\xmin_i+h_i)=\si(\xmin_i)$. Then,
$$\left\|\xmin_{S}+h_{S}\right\|_1-\left\|\xmin_{S}\right\|_1=h_S^\top\si(\xmin_S)$$
and
$$r(\xmin+h)-r(\xmin)\hspace{-0.05cm}=\hspace{-0.05cm}d\left\|h_{\Sc}\right\|_1+ h_S^\top[d\si(\xmin_S)-\xmin_S]-\frac{\left\|h\right\|_2^2}{2}.$$
%
Furthermore, since $Ah=A_Sh_S+A_{\Sc}h_{\Sc}$,
\begin{equation*}
\begin{split}
\fun(\xmin+h)-\fun(\xmin)&=\frac{1}{2}\left\|Ah\right\|_2^2+ h_{\Sc}^\top A_{\Sc}^\top(A\xmin-y)+\la d\left\|h_{\Sc}\right\|_1+\\
&+h_S^\top\left[\la d\si(\xmin_S)-\la\xmin_S+A_S^\top(A\xmin-y)\right]-\frac{\la}{2}\left\|h\right\|_2^2.
\end{split}
\end{equation*}
If $\si(\xmin_S)=\si(\xtrue_S)$, the quantity within the square parenthesis is the restricted subgradient which is null by construction in \eqref{defmin}.

%
Thus,
\begin{equation}\label{thenthen}
\begin{split}
\fun(\xmin+h)-\fun(\xmin)=\frac{1}{2}\left\|Ah\right\|_2^2-\frac{\la}{2}\left\|h\right\|_2^2+ h_{\Sc}^\top A_{\Sc}^\top(A\xmin-y)+\la d\left\|h_{\Sc}\right\|_1.
\end{split}
\end{equation}
We compute a lower bound for $h_{\Sc}^\top A_{\Sc}^\top(A\xmin-y)$, by exploiting the H\"{o}lder inequality:
\begin{equation}
\begin{split}
&|h_{\Sc}^\top A_{\Sc}^\top(A\xmin-y)|\leq \left\|h_{\Sc}\right\|_1\left\|A_{\Sc}^\top(A\xmin-y)\right\|_{\infty} \\
&\leq \left\|h_{\Sc}\right\|_1\left\|A_{\Sc}^\top A_S(\xmin_S-\xtrue_S)\right\|_{\infty}+\left\|h_{\Sc}\right\|_1\left\|A_{\Sc}^\top\noise\right\|_{\infty}.
\end{split}
\end{equation}
Since 
\begin{equation}\label{questa}
A_{\Sc}^\top A_S(\xmin_S-\xtrue_S)=A_{\Sc}^\top A_S(A_S^\top A_S)^{-1}(A_S^\top A_S)(\xmin_S-\xtrue_S),
\end{equation}
we elaborate on $A_S^\top A_S(\xmin_S-\xtrue_S)$ and exploit $(\omega,S)$-IRR to evaluate the term $\left\|A_{\Sc}^\top A_S(\xmin_S-\xtrue_S)\right\|_{\infty}$.
We have
\begin{equation}
 \begin{split}
 \|A_S^\top A_S(\xmin_S-\xtrue_S)\|_{\infty}&=\left\|A_S^\top A_S(A_S^\top A_S-\la I)^{-1}\left[\la \xtrue_S-\la d\si(\xtrue_S)+A_S^\top\noise\right]\right\|_{\infty}\\
 &\leq \left\|A_S^\top A_S(A_S^\top A_S-\la I)^{-1}\right\|_{\infty}\left[\la d(1-\mu)+\|A_S^\top\noise\|_{\infty}\right].
\end{split}
 \end{equation}
Hence, from \eqref{questa} and by using $(\omega,S)$-IRR, we get
\begin{equation*}
\begin{split}
&\left\|A_{\Sc}^\top A_S(\xmin_S-\xtrue_S)\right\|_{\infty}\leq\\
&~~~~~\leq\|A_{\Sc}^\top A_S^{\dagger\top}\|_\infty \left\|A_S^\top A_S(A_S^\top A_S-\la I)^{-1}\right\|_{\infty}\left[\la d(1-\mu)+\|A_S^\top\noise\|_{\infty}\right]\\
&~~~~~= \|A_S^{\dagger}A_{\Sc}\|_1 \left\|A_S^\top A_S(A_S^\top A_S-\la I)^{-1}\right\|_{\infty}\left[\la d(1-\mu)+\|A_S^\top\noise\|_{\infty}\right]\\
&~~~~~\leq \frac{1}{\omega} \left\|A_S^\top A_S(A_S^\top A_S-\la I)^{-1}\right\|_{\infty}\left[\la d(1-\mu)+\|A_S^\top\noise\|_{\infty}\right].
\end{split}
\end{equation*}
In conclusion, the following lower bound holds for \eqref{thenthen}:
\begin{equation*}
\begin{split}
&\fun(\xmin+h)-\fun(\xmin)\geq \frac{1}{2}\left\|Ah\right\|_2^2-\frac{\la}{2}\left\|h\right\|_2^2+q\left\|h_{\Sc}\right\|_1
\end{split}
\end{equation*}
where 
\begin{equation}
q=\la d-\frac{1}{\omega}\left\|A_S^\top A_S(A_S^\top A_S-\la I)^{-1}\right\|_{\infty}\left[\la d(1-\mu)+\|A_S^\top\noise\|_{\infty}\right]-\left\|A_{\Sc}^\top\noise\right\|_{\infty}.
\end{equation}

Given that $(\alpha,\phi,S)$-RE holds, with $\phi>\la$, we distinguish two cases.

If $h\in\cone(\alpha,S)\setminus\{0\}$, then $ \frac{1}{2}\left\|Ah\right\|_2^2-\frac{\la}{2}\left\|h\right\|_2^2\geq 0$. Thus, $\fun(\xmin+h)>\fun(\xmin)$ for any $h\in\cone(\alpha,S)\setminus\{0\}$ whenever $q>0$.

Otherwise, if $h\notin\cone(\alpha,S)\setminus\{0\}$, then $\|h_S\|_1\leq\frac{\|h_{\Sc}\|_1}{\alpha}$, which implies
\begin{equation*}
\begin{split}
\|h\|_2^2&\leq \varepsilon \|h\|_1 = \varepsilon \|h_S\|_1 +\varepsilon \|h_{\Sc}\|_1\leq \varepsilon \|h_{\Sc}\|_1 \frac{\alpha+1}{\alpha}.
\end{split}
\end{equation*}
Therefore, $\fun(\xmin+h)-\fun(\xmin)\geq \left\|h_{\Sc}\right\|_1\left( q-\la\varepsilon\frac{1+\alpha}{\alpha}\right)$.

As $\|h_S\|_1\leq\frac{\|h_{\Sc}\|_1}{\alpha}$, if $h\neq 0$, then $h_{\Sc}\neq 0$.
Thus, if $h\neq 0$, then $\fun(\xmin+h)>\fun(\xmin)$  whenever  C3 holds.
 \end{proof}

This result yields some considerations. 
\begin{remark}\label{rem:relax}
In the noise-free case $\eta=0$, we can consider $\la$ arbitrarily small. Then, in condition C3 becomes $d-\frac{1}{\omega} d(1-\mu)>\varepsilon\frac{\alpha+1}{\alpha}$. Given that $\varepsilon$ is arbitrarily small, we must have $\omega>1-\mu$. This irrepresentable condition is weaker with respect to Lasso, which requires $\omega\geq 1$. In particular, a larger $\mu$ implies a more relaxed irrepresentable condition.

In other terms, the irrepresentable condition for MCPS$^2$ is tuned based on the minimum non-zero magnitude in $\xtrue$. Interestingly, in the limit case  $\mu=1$, no irrepresentable condition is required; see Section \ref{sec:ta3} for  details.

Instead, in the presence of noise the value of $\la$ must be sufficiently large to keep $q>0$. On the other hand, condition \eqref{limit_lambda}
$$\|(A_S^\top A_S-\la I)^{-1}\|_{\infty}\left[ \la d (1-\mu)+\|A_S^\top\noise\|_{\infty}\right]\leq  \min_{i\in S}|\xtrue_i|$$
imposes an upper bound on $\la$ depending on $\max_{i\in S}|\xtrue_i|$, which limits the resilience to noise.
This is similar to condition \eqref{lasso_limitazione_su_lambda} for Lasso $$ \|(A_S^\top A_S)^{-1}\|_{\infty}(\|A_S^\top \noise\|_{\infty}+\la)\leq \min_{i\in S}|\xtrue_i|$$ which limits the value of $\la$ baed on the noise and $\max_{i\in S}|\xtrue_i|$ as well. However, condition \eqref{limit_lambda} may be weaker than \eqref{lasso_limitazione_su_lambda} when $\mu$ is close to $1$.
\end{remark}

 To assess the considerations in Remark \eqref{rem:relax}, we present some numerical examples. We consider some matrices $A$ and we evaluate their irrepresentable properties in terms of Lasso, see \eqref{lasso_irr}, and MCPS$^2$, see condition C3 in Theorem \ref{theo1}.
 We assume that  $\xtrue$ has length $n=100$ and sparsity $k=5$ and we randomly generate the support. The non-zero entries of $\xtrue$ have random magnitudes in $[\frac{1}{2},1]$, that is, $d=1$ and $\mu=\frac{1}{2}$. We consider $A\in\R^{m,n}$ with Gaussian independent entries $\mathcal{N}(0,\frac{1}{m})$, $m\in[10,100]$. We illustrate both the noise-free case and a bounded measurement noise $\noise\in\R^m$ with $\|\noise\|_{\infty}=10^{-3}$. For each $m$, we generate 1000 matrices $A$ and we check the irrepresentable condition rate, that is, how many times the irrepresentable conditions \eqref{lasso_irr} and C3 hold, given a suitable tuning of $\la$, guaranteeing the VSC.
\begin{figure}[ht!]
	\centering
	\includegraphics[width=0.6\columnwidth]{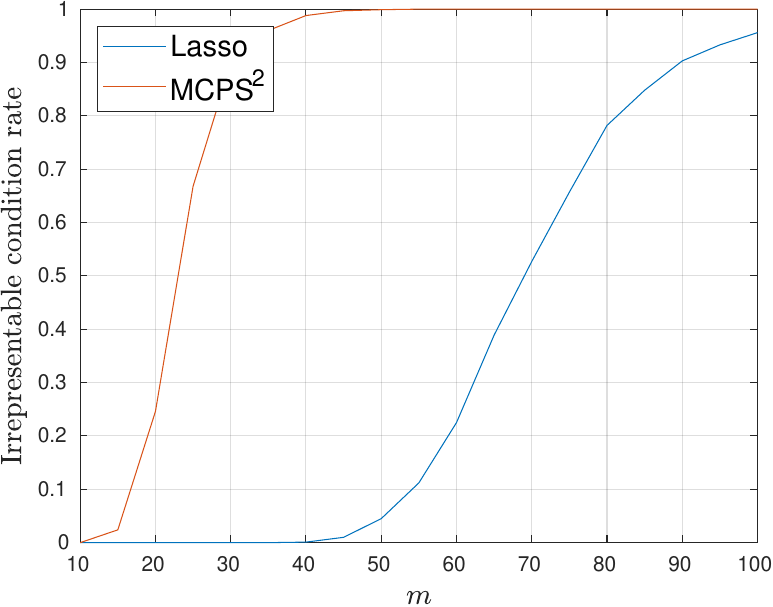}\qquad
	\includegraphics[width=0.6\columnwidth]{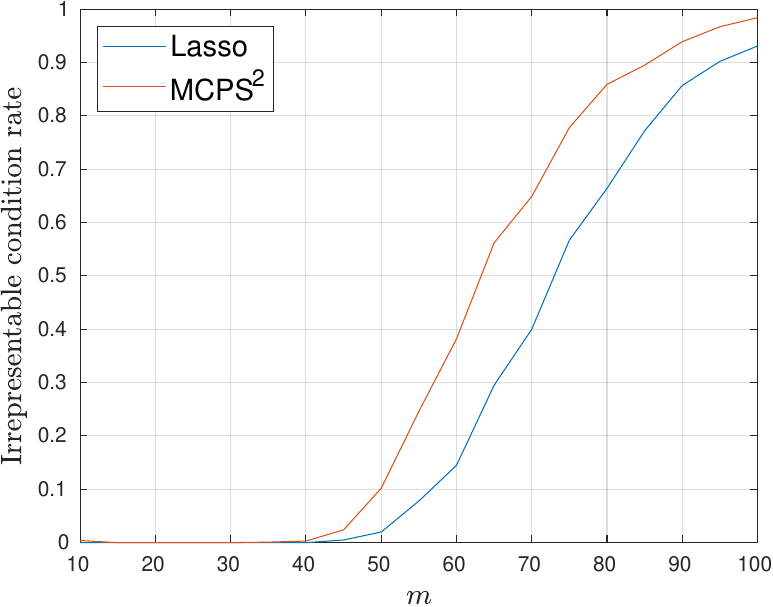}
  \vskip0.5cm
	\caption{Irrepresentable condition rate \eqref{lasso_irr} (Lasso) and C3 (MCPS$^2$), over 1000 runs, for $\noise=0$ (top) and $\|\noise\|_{\infty}=10^{-3}$ (bottom). In both cases, the irrepresentable condition that guarantees the VSC of MCPS$^2$ is satisfied more often than that of Lasso, which encourage the use of MCPS$^2$.}
	\label{fig:0}
\end{figure}
In Fig. \ref{fig:0}, we show the results of this test, which confirm that the irrepresentable condition that guarantees the VSC of MCPS$^2$ is satisfied more frequently than that of Lasso. This is particularly evident in the noise-free case, where 40 measurements are sufficient to estimate the correct support, provided that we run an algorithm that finds the corresponding minimum of MCPS$^2$.

%

\begin{remark}
Theorem \ref{theo1} states the possibility of achieving the desired $\xmin$ defined in \eqref{defmin} in the sense that $\xmin$ is a local minimum of the proposed functional, therefore we can reach it by a descent algorithm, given a suitable starting point.
\end{remark}

\begin{remark}
The value of $\alpha$ is not assessed. Actually, since $\varepsilon$ is arbitrarily small, in C3 we could neglect the term $\varepsilon(\alpha+1)$, thus remove $\alpha$. For Lasso, ($\alpha,\phi,S$)-RE is useful when $\alpha=3$, see \cite[Theorem 11.1]{has15book}. Moreover, the bound on the $\ell_2$ error in \cite[Theorem 11.1]{has15book} is controlled by $\frac{3\la\sqrt{k}}{\phi}$, which suggests that $\phi\gg\la$ is necessary to have a small error. Instead, to prove Theorem \ref{theo1}, we only require $\phi\geq\la$. 
\end{remark}
%
%
%

\subsection{Global minimum}\label{sec:ta2}
In this section, we provide sufficient conditions such that $\xmin$ defined in Theorem \ref{theo1} is not only a local minimum, but also the global minimum of $\fun$.
From Theorem \ref{theo1}, let us assume that $\xmin$ is the unique minimum over  $\ball:=\{\xmin+h, \text{ for all } h\in\R^n \text{ with }\|h\|_{\infty}\leq \varepsilon\}$ for some $\varepsilon>0$. Then, if the global minimum lies in $\ball$, it necessarily corresponds to $\xmin$. In the following proposition, we provide sufficient conditions for this occurrence.
\begin{proposition}\label{prop1}
Let $\theta=\|(y-A\xmin)^\top A\|_{\infty}$. Let us assume $\theta<\la d$. If $A$ satisfies $(\alpha,\phi,S)$-RE with $\alpha\geq \frac{2\la d+\theta}{\la d-\theta}$, $\phi> \la$, and  $\frac{2(\theta+2\la d)\sqrt{k}}{\phi-\la}\leq\varepsilon$, then, $\xmin$ defined in \eqref{defmin} is the global minimum of $\fun$.
\end{proposition}
\begin{proof}
 Let $\xgmin$ be the  global minimum of $\fun$, and $\nu:=\xgmin-\xmin$. Then, $\fun(\xmin+\nu)=\fun(\xgmin)\leq \fun(\xmin)$.
 This is equivalent to
%
$\frac{1}{2}\|A\nu+A\xmin-y\|_2^2+ \la d \|\xmin+\nu\|_1-\frac{\la}{2}\|\xmin+\nu\|_2^2\leq \frac{1}{2}\|A\xmin-y\|_2^2+ \la d \|\xmin\|_1-\frac{\la}{2}\|\xmin\|_2^2$.
Further simplifications lead to 
 $\frac{1}{2}\|A\nu\|_2^2+(A\xmin-y)^\top A\nu\leq \la d \|\xmin\|_1-\la d \|\xmin+\nu\|_1+\frac{\la}{2}\|\nu\|_2^2+\la \nu_S^\top\xmin_S$.
 Now, we remark that
 \begin{equation*}
  \begin{split}
   \|\xmin\|_1-\|\xmin+\nu\|_1& = \|\xmin_S\|_1-\|\xmin_S+\nu_S\|_1 -\|\nu_{\Sc}\|_1\\
  &\leq \|\nu_S\|_1 -\|\nu_{\Sc}\|_1.
  \end{split}
 \end{equation*}
 Moreover, since $\xmin_S=\xgmin_S-\nu_S$, and $\|\xgmin\|_{\infty}\leq d$,
 \begin{equation*}
  \begin{split}
  \|\nu\|_2^2+2\nu_S^\top\xmin_S&= \|\nu_S\|_2^2+\|\nu_{\Sc}\|_2^2+2\nu_S^\top\xgmin_S -2\|\nu_S\|_2^2\\
&\leq -\|\nu_S\|_2^2+\|\nu_{\Sc}\|_2^2+2d\|\nu_S\|_1.
   \end{split}
 \end{equation*}
 Then,
 \begin{equation*}
 \frac{1}{2}\|A\nu\|_2^2\leq (\theta+\la d) \|\nu_S\|_1+\theta\|\nu_{\Sc}\|_1+\la r( v_S)-\la r(v_{\Sc}).
\end{equation*}
If $A$ satisfies $(\alpha,\phi,S)$-RE, then
$\frac{\phi+\la}{2}\|\nu_S\|_2^2+\frac{\phi-\la}{2}\|\nu_{\Sc}\|_2^2 \leq(\theta+2\la d) \|\nu_S\|_1+(\theta-\la d)\|\nu_{\Sc}\|_1.$
If $\phi>\la$, then $0\leq\frac{\phi+\la}{2}\|\nu_S\|_2^2+\frac{\phi-\la}{2}\|\nu_{\Sc}\|_2^2$, and
$(\la d-\theta)\|\nu_{\Sc}\|_1\leq (\theta+2\la d)\|\nu_{S}\|_1.$
Therefore, if $\theta <\la d$, then $\nu\in \cone\left(\frac{2\la d+\theta}{\la d-\theta}, S\right)$.
 Moreover,
\begin{equation}\label{final_estimate_global}
\begin{split}
\frac{\phi-\la}{2}\|\nu\|_2^2&\leq (\theta+2\la d)\|\nu_S\|_1\leq (\theta+2\la d)\sqrt{k}\|\nu\|_2\\
\Rightarrow \|\nu\|_2&\leq \frac{2(\theta+2\la d)\sqrt{k}}{\phi-\la}.
 \end{split}
\end{equation}
If $\|\nu\|_2\leq \varepsilon$, then $\xgmin\in \ball$. Since $\xgmin$ is the global minimum, while $\xmin$ is the global minimum limited to $\ball$, then $\|\nu\|_2\leq \varepsilon$ implies $\xgmin=\xmin$.
According to \eqref{final_estimate_global}, a sufficient condition to have $\|\nu\|_2\leq \varepsilon$ is
$\frac{2(\theta+2\la d)\sqrt{k}}{\phi-\la}\leq\varepsilon.$
\end{proof}
The bound in \eqref{final_estimate_global} might be refined with further computations, which will be proposed in a future extended work. However, even though not perfectly tight, it well illustrates the robustness to noise: if $\noise$ decreases, and as a consequence $\la$ and $\theta$ can be proportionally decreased, then also $\|\nu\|_2$ decreases, which increases the probability that $\xmin$ is the global minimum.

\subsection{Noise-free, binary case}\label{sec:ta3}
To conclude, we illustrate the noise-free case  $\noise=0$ with $\mu=1$, i.e., $\xtrue_S\in\{-d,0,d\}^k$. From Theorem \ref{theo1} and Proposition \ref{prop1}, we derive the following result.
\begin{corollary}\label{cor1}
Let us assume $\noise=0$ and $\mu=1$.
If  $A$ satisfies $(\omega,S)$-IRR, and $A$ satisfies $(\alpha,\phi,S)$-RE with $\phi>\la$, and $\alpha=2$, then
 $\xmin=\xtrue$ is a local minimum of $\fun$. More precisely, $\xmin=\xtrue$ is the unique minimum in $\ball$, with $\varepsilon=\frac{2}{3}d\omega$.
Moreover, if $\la <\frac{\omega \phi}{6\sqrt{k}+\omega}$, then $\xtrue$ is the global minimum of $\fun$.
\end{corollary}
The proof can be straightforwardly obtained by replacing $\noise=0$ and $\mu=1$ in Theorem \ref{theo1} and Proposition \ref{prop1}. Corollary \ref{cor1} shows that in the favorable case without noise and with extreme non-zero values, the proposed approach is always effective, since it is sufficient to set a sufficiently small $\la$ to obtain that the true $\xtrue$ is the global minimum of $\fun$, without bias. For more details on CS with discrete-valued signals, see \cite{fox18tsp}.
\subsection{Discussion}
The proposed analysis is theoretical, as some parameters, such as $\omega$ and $\phi$, are not a priori known, and therefore the choice of the design parameter $\la$ is not precisely determined. The main contribution of this analysis is to state that MCPS$^2$, differently from Lasso, does not require the classical IRR. Similarly to classical CS, in future work, conditions on $A$ will be studied that are provable in practice, at least for matrices' ensembles, such as the restricted isometry property. From a practical viewpoint, we expect that the relaxation of IRR leads to a reduction of the number $m$ of measurements needed for a perfect recovery of the support; this is illustrated in numerical simulations in the next section.

Finally, we remark that finding the global minimum of $\fun$ may be not straightforward, since $\fun$ is non-convex. In the next section, we test different algorithms to achieve the global minimum, and show their effectiveness through numerical simulations.
\section{Algorithms and numerical simulations}\label{sec:nr}
In this section, we test different algorithms to achieve the minimum of $\fun$ as defined in \eqref{mcps2}; this task is challenging due to non-convexity. We propose two approaches: the semidefinite programming relaxation (SDR), supported by recent results on polynomial optimization \cite{lasbook}, and the alternating direction method of multipliers (ADMM, \cite{boy10}).
\subsection{Semi-algebraic optimization}\label{sec:POP}
Since $\fun$ in \eqref{mcps2} is semi-algebraic, the theory developed in \cite{las01} can be applied to compute the global minimum. In a nutshell, given a polynomial or semi-algebraic optimization problem, a hierarchy  of SDR's can be constructed, whose solutions converge to the global optimal solution. The hierarchy generically has finite convergence, see, e.g., \cite{nie14}. Therefore, the global minimum can be achieved by solving an SDR of sufficiently large order.
A shortcoming of the SDR approach is the numerical complexity, which is of order $\mathcal{O}(n^{\zeta})$, $n$ being the number of variables and $\zeta$ the relaxation order. For this motivation, in this paper, we consider only the SDR of order $\zeta=1$, which corresponds to the Shor's relaxation.
\subsection{ADMM}\label{sec:MADMM}
ADMM is an iterative algorithm, widely used in convex optimization for its fast  convergence and simplicity of implementation, see \cite{boy10}. In the non-convex setting, the convergence of ADMM to a local minimum has been proven only for some classes of functionals. In particular, in \cite{hon16}, convergence is proven for non-convex functionals that can be split into the sum of a non-convex, smooth term and of a convex, not necessarily smooth, term. $\fun$ in \eqref{mcps2} has this property, i.e., we can write it as
\begin{equation*}
\begin{split}
&\min_{x,z\in \D} \frac{1}{2}\left\|y-A x \right\|_2^2-\frac{\la}{2}\left\|x \right\|_2^2+\la d \left\|z\right\|_1~\text{ s. t. }~ z=x
\end{split}
\end{equation*}
with
The associated augmented Lagrangian is:
$\lag(x,z)=\frac{1}{2}\left\|y-A x \right\|_2^2-\frac{\la}{2}\left\|x \right\|_2^2+\la d \left\|z\right\|_1+u^\top(x-z)+\frac{\rho}{2}\left\|x-z \right\|_2^2$
where $u\in\R^n$ is the dual variable, and $\rho>0$. Then, we apply ADMM as explained in \cite[Section 2]{hon16}, which consists in iteratively minimizing $\lag$ with respect to $x$ and to $z$, and updating $u$. This procedure is summarized in Algorithm \ref{alg:MADMM}, where $P$ denotes the  operator that projects onto $[-d, d]^n$, and $\soft_a$ is the soft thresholding operator.
\begin{algorithm}
\setstretch{1.2}
     \renewcommand{\algorithmicrequire}{\textbf{Input:}}
    \renewcommand{\algorithmicensure}{\textbf{Output:}}
  \caption{ADMM for MCPS$^2$}\label{alg:MADMM}
  \begin{algorithmic}[1] 
   \REQUIRE $A,y=A\xtrue+\noise$, $\la>0$, $\rho>0$
 \ENSURE  $x_{T_{stop}}$ = estimate of $\xtrue$\\

    \STATE Initialize $z_0,w_0\in\R^n$
    \FORALL{$t=1,\dots,T_{stop}$}
    \STATE  $x_t=\argmin{x\in\R^n} \lag(x,z_{t-1})$\\
     $=\left[A^\top A + (\rho-\la)I\right]^{-1}\big( A^\top y + \rho z_{t-1} - u_{t-1}\big)$
    %
    \STATE $z_t=\argmin{z\in[-d,d]^n} \lag(x_t,z)=P\left(\soft_{\frac{\la d}{\rho}}\big(x_t+\frac{u_{t-1}}{\rho}\big)\right)$
    %
    \STATE $u_t=u_{t-1}+\rho(x_t-z_t)$
    \ENDFOR     
  \end{algorithmic}
\end{algorithm} 
%

\subsection{Numerical results}
\begin{figure}[ht!]
	\centering
	\includegraphics[width=0.66\columnwidth]{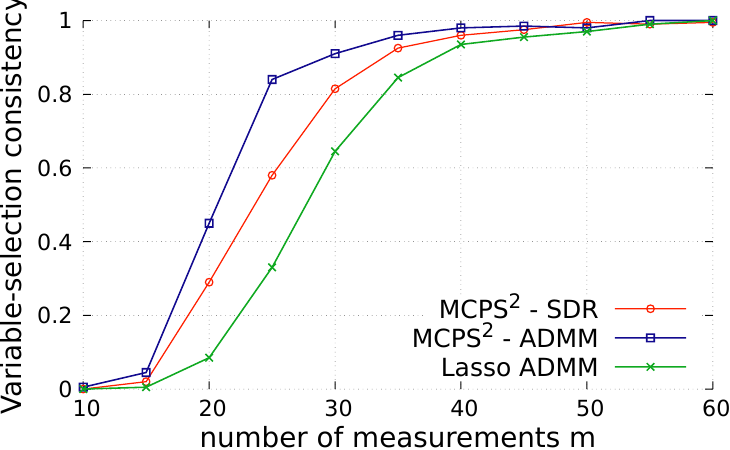}
	\vskip0.6cm
	\includegraphics[width=0.66\columnwidth]{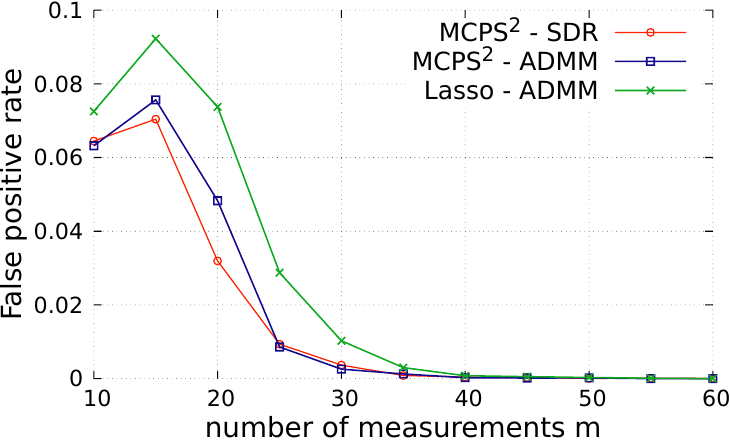}
	\vskip0.6cm
	\includegraphics[width=0.66\columnwidth]{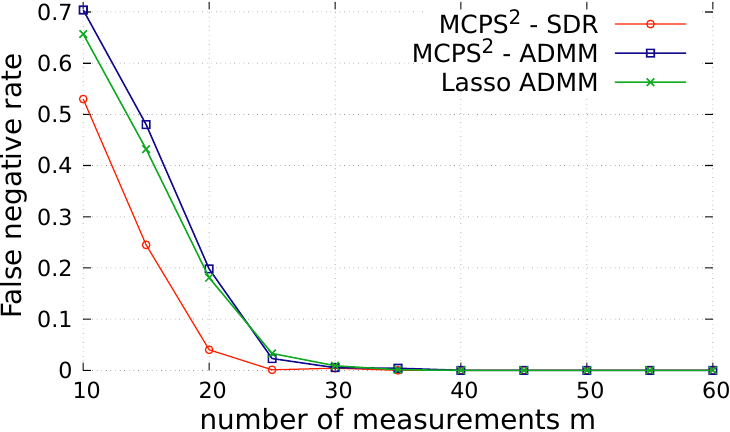}
	\caption{Comparison between the proposed MCPS$^2$ optimization problem, solved via SDR and ADMM, and the Lasso, solved via ADMM in terms of VSC}
	\label{fig:1}
\end{figure}

\begin{figure}[ht!]
\centering
\includegraphics[width=0.66\columnwidth]{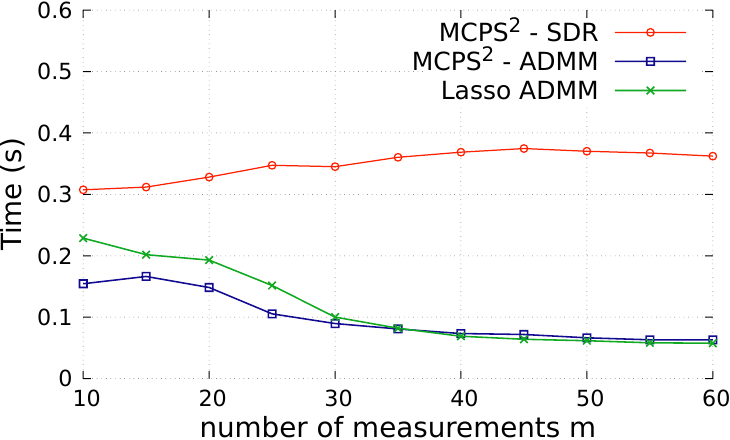}
	\caption{Comparison between the proposed MCPS$^2$ optimization problem, solved via SDR and ADMM, and the Lasso, solved via ADMM: convergence time}
	\label{fig:2}
\end{figure}
We illustrate some numerical simulations with SDR and ADMM approaches to solve MCPS$^2$. In particular, we compare them to Lasso, solved with ADMM.

The considered approaches are compared in terms of VSC, i.e., the number of experiments where the support is exactly recovered, and false positive/negative rate, that is the rate of zeros estimated as non-zeros, and vice-versa. Finally, we compare the runtimes. All the algorithms are implemented in C++; for SDR, we use the Mosek C++ Fusion API, see \cite{mosek}.

The results, averaged over 200 runs, are shown in Fig. \ref{fig:1}. We see that MCPS$^2$ performs better than Lasso in terms of VSC, as expected from the proposed theoretical results in Sec. \ref{sec:ta}. The enhancement is particularly evident for $m\in[20,35]$. For example, for $m=30$, Lasso is not reliable, with the $60\%$ of successful support recovery, while MCPS$^2$ attains the $90\%$. Furthermore, we can see that MCPS$^2$ outperforms Lasso in terms of false positive rate, while the false negative rate of SDR is better than that of ADMM.

In Fig. \ref{fig:2}, we show the convergence time (in seconds), which is sufficiently fast for all the considered methods. As expected, SDR is slower than ADMM.
%

As to MCPS$^2$, SDR and ADMM algorithms both suboptimal: on the one hand, SDR can be improved by increasing the relaxation order, now set to 1 to minimize the runtime; on the other hand, for ADMM, suitable initial conditions could be investigated to achieve the desired minimum. Moreover,in the proposed experiments, we observe that when $k$ is correctly estimated, then the achieved point corresponds to $\xmin$. This suggests that $\xmin$ might be the unique minimum with sparsity $k$; in case, given the knowledge of $k$, this could be used to verify whether the achieved point is $\xmin$, and, if not, ADMM should be run again with different initial conditions to search $\xmin$. These points will be studied in future extended work.

%
\section{Conclusions}\label{sec:con}
In this work, we analyse MCPS$^2$, a non-convex, semi-algebraic optimization problem for sparse learning. Specifically, we study its variable-selection consistency in a compressed sensing framework, and  we prove that, differently from  Lasso, MCPS$^2$ relaxes the irrepresentable condition. In practice, this implies that MCPS$^2$ requires less measurements than  Lasso. Future work will be oriented to provide conditions guaranteeing the variable-selection consistency of MCPS$^2$ that can be a priori verified. and to develop optimal strategies to achieve the global minimum.

\bibliographystyle{IEEEtran}
\bibliography{refs}

\begin{thebibliography}{10}
\providecommand{\url}[1]{#1}
\csname url@samestyle\endcsname
\providecommand{\newblock}{\relax}
\providecommand{\bibinfo}[2]{#2}
\providecommand{\BIBentrySTDinterwordspacing}{\spaceskip=0pt\relax}
\providecommand{\BIBentryALTinterwordstretchfactor}{4}
\providecommand{\BIBentryALTinterwordspacing}{\spaceskip=\fontdimen2\font plus
\BIBentryALTinterwordstretchfactor\fontdimen3\font minus
  \fontdimen4\font\relax}
\providecommand{\BIBforeignlanguage}[2]{{%
\expandafter\ifx\csname l@#1\endcsname\relax
\typeout{** WARNING: IEEEtran.bst: No hyphenation pattern has been}%
\typeout{** loaded for the language `#1'. Using the pattern for}%
\typeout{** the default language instead.}%
\else
\language=\csname l@#1\endcsname
\fi
#2}}
\providecommand{\BIBdecl}{\relax}
\BIBdecl

\bibitem{has09book}
T.~Hastie, R.~Tibshirani, and J.~Friedman, \emph{The Elements of Statistical
  Learning}, 2nd~ed.\hskip 1em plus 0.5em minus 0.4em\relax Springer New York
  Inc., 2009.

\bibitem{has15book}
T.~Hastie, R.~Tibshirani, and M.~Wainwright, \emph{Statistical Learning with
  Sparsity: The {L}asso and Generalizations}, 2nd~ed.\hskip 1em plus 0.5em
  minus 0.4em\relax CRC press, 2015.

\bibitem{bru19book}
S.~L. Brunton and J.~N. Kutz, \emph{Data-Driven Science and Engineering:
  Machine Learning, Dynamical Systems, and Control}.\hskip 1em plus 0.5em minus
  0.4em\relax Cambridge University Press, 2019.

\bibitem{car13}
A.~Y. Carmi, ``Compressive system identification: Sequential methods and
  entropy bounds,'' \emph{Dig. Signal Process.}, vol.~23, no.~3, pp. 751--770,
  2013.

\bibitem{fat18}
S.~{Fattahi} and S.~{Sojoudi}, ``Data-driven sparse system identification,'' in
  \emph{Proc. Allerton Conf. Commun. Control Comput.}, 2018, pp. 462--469.

\bibitem{bru16}
S.~L. Brunton, J.~L. Proctor, and J.~N. Kutz, ``Discovering governing equations
  from data by sparse identification of nonlinear dynamical systems,''
  \emph{PNAS}, vol. 113, no.~15, pp. 3932--3937, 2016.

\bibitem{gal16book}
M.~Gallieri, \emph{{Lasso-MPC} - Predictive Control with $\ell_1$-Regularised
  Least Squares}.\hskip 1em plus 0.5em minus 0.4em\relax Springer Int., 2016.

\bibitem{tar18}
E.~Tartaglione, S.~Leps\o{}y, A.~Fiandrotti, and G.~Francini, ``Learning sparse
  neural networks via sensitivity-driven regularization,'' in \emph{Proc. Int.
  Conf. Neural Inf. Process. Sys. (NIPS)}, 2018, pp. 3882--3892.

\bibitem{bol19}
H.~B{\"o}lcskei, P.~Grohs, G.~Kutyniok, and P.~Petersen, ``Optimal
  approximation with sparsely connected deep neural networks,'' \emph{SIAM J.
  Math. Data Sci.}, vol.~1, no.~1, pp. 8--45, 2019.

\bibitem{don06}
D.~L. Donoho, ``Compressed sensing,'' \emph{IEEE Trans. Inf. Theory}, vol.~52,
  no.~4, pp. 1289--1306, 2006.

\bibitem{fou13}
S.~Foucart and H.~Rauhut, \emph{A Mathematical Introduction to Compressive
  Sensing}.\hskip 1em plus 0.5em minus 0.4em\relax New York: Springer, 2013.

\bibitem{tot11}
R.~T\'{o}th, B.~M. Sanandaji, K.~Poolla, and T.~L. Vincent, ``Compressive
  system identification in the linear time-invariant framework,'' in
  \emph{Proc. IEEE Conf. Decis. Control (CDC)}, 2011, pp. 783--790.

\bibitem{san11}
B.~M. Sanandaji, T.~L. Vincent, M.~B. Wakin, and R.~T\'{o}th, ``Compressive
  system identification of lti and ltv arx models,'' in \emph{Proc. IEEE Conf.
  Decis. Control (CDC)}, 2011, pp. 783--790.

\bibitem{tib96}
R.~Tibshirani, ``Regression shrinkage and selection via the lasso,'' \emph{J.
  Roy. Stat. Soc. Series B}, vol.~58, pp. 267--288, 1996.

\bibitem{fuc04}
J.~J. Fuchs, ``On sparse representations in arbitrary redundant bases,''
  \emph{IEEE Trans. Inf. Theory}, vol.~50, no.~6, pp. 1341--1344, 2004.

\bibitem{zha06}
P.~Zhao and B.~Yu, ``On model selection consistency of {L}asso,'' \emph{J.
  Mach. Learn. Res.}, vol.~7, pp. 2541 -- 2563, 2006.

\bibitem{wai09}
M.~J. Wainwright, ``Sharp thresholds for high-dimensional and noisy sparsity
  recovery using {$\ell_1$}-constrained quadratic programming ({L}asso),''
  \emph{IEEE Trans. Inf. Theory}, vol.~55, no.~5, pp. 2183--2202, 2009.

\bibitem{bul11}
P.~B{\"u}hlmann and S.~van~de Geer, \emph{Variable selection with the
  Lasso}.\hskip 1em plus 0.5em minus 0.4em\relax Springer Berlin Heidelberg,
  2011, pp. 183--247.

\bibitem{can08rew}
E.~J. Cand\`es, M.~B. Wakin, and S.~Boyd, ``Enhancing sparsity by reweighted
  $\ell_1$ minimization,'' \emph{J. Fourier Anal. Appl.}, vol.~14, no. 5-6, pp.
  877--905, 2008.

\bibitem{woo16}
J.~Woodworth and R.~Chartrand, ``Compressed sensing recovery via nonconvex
  shrinkage penalties,'' \emph{Inverse Problems}, vol.~32, no.~7, pp.
  75\,004--75\,028, 2016.

\bibitem{sel17}
I.~Selesnick, ``Sparse regularization via convex analysis,'' \emph{IEEE Trans.
  Signal Process.}, vol.~65, no.~17, pp. 4481--4494, 2017.

\bibitem{fox19cdc}
V.~{Cerone}, S.~M. {Fosson}, and D.~{Regruto}, ``Sparse linear regression with
  compressed and low-precision data via concave quadratic programming,'' in
  \emph{Proc. Conf. Decis. Control (CDC)}, 2019, pp. 6971--6976.

\bibitem{lasbook}
J.-B. Lasserre, \emph{An introduction to polynomial and semi-algebraic
  optimization}.\hskip 1em plus 0.5em minus 0.4em\relax Cambridge University
  Press, 2015.

\bibitem{tib13}
R.~J. Tibshirani, ``{The {L}asso problem and uniqueness},'' \emph{Electronic
  Journal of Statistics}, vol.~7, pp. 1456--1490, 2013.

\bibitem{fuc05}
J.~J. Fuchs, ``Recovery of exact sparse representations in the presence of
  bounded noise,'' \emph{IEEE Trans. Inf. Theory}, vol.~51, no.~10, pp.
  3601--3608, 2005.

\bibitem{mei06}
N.~Meinshausen and P.~B\"{u}hlmann, ``High-dimensional graphs and variable
  selection with the {L}asso,'' \emph{Ann. Stat.}, vol.~34, pp. 1436--1462,
  2006.

\bibitem{fou09}
S.~Foucart and M.-J. Laui, ``Sparsest solutions of underdetermined linear
  systems via $\ell_q$ minimization for $0<q\leq 1$,'' \emph{Appl. Comput.
  Harmon. Anal.}, vol.~26, pp. 395--407, 2009.

\bibitem{zha10MCP}
C.-H. Zhang, ``Nearly unbiased variable selection under minimax concave
  penalty,'' \emph{Ann. Statist.}, vol.~38, no.~2, pp. 894--942, 2010.

\bibitem{fox18}
S.~M. Fosson, ``A biconvex analysis for lasso $\ell_1$ reweighting,''
  \emph{IEEE Signal Process. Lett.}, vol. early access, no.~nn, pp. 1 -- 1,
  2018.

\bibitem{fox18tsp}
------, ``Non-convex {L}asso-kind approach to compressed sensing for
  finite-valued signals,'' \emph{arxiv.org/abs/1811.03864v2}, 2018.

\bibitem{van09}
S.~A. van~de Geer and P.~B\"{u}hlmann, ``On the conditions used to prove oracle
  results for the {L}asso,'' \emph{Electron. J. Stat.}, vol.~3, pp. 1360--1392,
  2009.

\bibitem{ras10}
G.~Raskutti, M.~J. Wainwright, and B.~Yu, ``Restricted eigenvalue properties
  for correlated {G}aussian designs,'' \emph{J. Machine Learn. Res.}, vol.~11,
  pp. 22\,641--2259, 2010.

\bibitem{boy10}
S.~Boyd, N.~Parikh, E.~Chu, B.~Peleato, and J.~Eckstein, ``Distributed
  optimization and statistical learning via the alternating direction method of
  multipliers,'' \emph{Found. Trends Mach. Learn.}, vol.~3, no.~1, pp. 1 --
  122, 2010.

\bibitem{las01}
J.~Lasserre, ``Global optimization with polynomials and the problem of
  moments,'' \emph{SIAM J. Optim.}, vol.~11, no.~3, pp. 796--817, 2001.

\bibitem{nie14}
J.~Nie, ``Optimality conditions and finite convergence of lasserre's
  hierarchy,'' \emph{Math. Program.}, vol. 146, no.~1, pp. 97--121, 2014.

\bibitem{hon16}
M.~Hong, Z.~Q. Luo, and M.~Razaviyayn, ``Convergence analysis of alternating
  direction method of multipliers for a family of nonconvex problems,''
  \emph{SIAM J. Optim.}, vol.~26, no.~1, pp. 337--364, 2016.

\bibitem{mosek}
\BIBentryALTinterwordspacing
{MOSEK}, \emph{{F}usion {API} for {C}++, version 9.1.5}, 2019. [Online].
  Available: \url{www.mosek.com/documentation/}
\BIBentrySTDinterwordspacing

\end{thebibliography}

\end{document}